\documentclass[12pt, a4paper]{article}
\usepackage[left=1in,right=1in,top=1.0in,bottom=1in]{geometry}

\usepackage{amsmath,amssymb,amsthm,amsfonts,afterpage,hyperref}
\usepackage{amscd,subeqnarray}

\usepackage{algorithm}
\usepackage{graphicx}
\usepackage{authblk}
\usepackage{tikz}
\usepackage{array}
\usepackage{wrapfig}
\usepackage{color}
\usepackage{ragged2e}
\usepackage{stmaryrd}
\usepackage{siunitx}
\usepackage{multirow}
\usepackage{xcolor,cancel}
\usepackage{boldfonts}
\usepackage{symbols}
\usepackage{footmisc}
\usepackage{float}
\usepackage{setspace}
\usepackage{bm}
\usepackage{booktabs}
\usepackage{tikz, xcolor}
\usepackage{soul}
\usepackage{algpseudocode}

\usepackage[sort, numbers, compress]{natbib}

\usetikzlibrary{shapes,arrows, positioning,calc}	
\usetikzlibrary{arrows,snakes,backgrounds}

\usepackage{enumitem}

\newlist{notes}{enumerate}{1}
\setlist[notes]{label=Note: ,leftmargin=*}

\usepackage[english]{babel}	
\usepackage{import}
\usepackage{caption}
\usepackage{subcaption}
\usepackage{mathtools}
\usepackage{enumitem}
\usetikzlibrary{automata, positioning}
\usepackage{fancyhdr}
\usepackage{url}
\usepackage[mathscr]{eucal}
\usepackage{acro}

\DeclareMathAlphabet{\mathpzc}{OT1}{pzc}{m}{it}

\newtheorem{definition}{Definition}

\providecommand{\keywords}[1]
{
  \small	
  \textbf{\textit{Keywords---}} #1
}

\title{Neural solver for sixth-order ordinary differential equations  }

\author[1,2]{Janavi Bhalala}
\author[3]{B. Veena S. N. Rao \thanks{Faculty Advisor}}

\affil[1]{Department of Computer Science, The University of Texas at Dallas, Richardson, TX-75080, USA}
\affil[2]{Department of Computer Science, Texas A\&M University-Corpus Christi, Corpus Christi, TX-78412, USA}
\affil[3]{Department of Mathematics \& Statistics,
Texas A\&M University - Corpus Christi
Corpus Christi, TX 78412, USA}
\date{}

\begin{document}

\maketitle

\begin{abstract}
A method for approximating sixth-order ordinary differential equations is proposed, which utilizes a deep learning feedforward artificial neural network, referred to as a neural solver. The efficacy of this unsupervised machine learning method is demonstrated through the solution of two distinct boundary value problems (BVPs), with the method being extended to include the solution of a sixth-order ordinary differential equation (ODE). The proposed mean squared loss function is comprised of two terms: the differential equation is satisfied by the first term, while the initial or boundary conditions are satisfied by the second. The total loss function is minimized using a quasi-Newton optimization method to obtain the desired network output. The approximation capability of the proposed method is verified for sixth-order ODEs. Point-wise comparisons of the approximations show strong agreement with available exact solutions. The proposed algorithm minimizes the overall learnable network hyperparameters in a given BVP.  Simple minimization of the total loss function yields highly accurate results even with a low number of epochs. Therefore, the proposed framework offers an attractive setting for the computational mathematics community.
\end{abstract}

\noindent \keywords{Machine learning, Artificial neural networks,  sixth-order ordinary differential equations, Optimizer, Activation Function }


\section{Introduction}
The study of higher-order differential equations is a critical pursuit in numerous scientific and engineering fields, as these equations serve as the mathematical language for describing complex physical phenomena \cite{hyun_mms_2021,yoon2022finite,gou2015modeling,Mallikarjunaiah2015,ferguson2015numerical}.  While analytical solutions are invaluable for their theoretical elegance, they are often nonexistent or exceedingly difficult to obtain for the nonlinear and high-dimensional problems that frequently arise in real-world applications. This intractability necessitates a reliance on numerical methods. The development of robust and accurate numerical solvers for higher-order differential equations is therefore a fundamental and pressing challenge. These methods provide a powerful computational framework for approximating solutions to problems that lack closed-form expressions, thereby enabling detailed quantitative analysis and prediction. The ability to simulate and predict the behavior of systems, ranging from the intricate quantum mechanical interactions to the complex dynamics of fluid flow in biological and mechanical systems, hinges on our capacity to find high-fidelity numerical solutions. The pursuit of novel and more efficient numerical techniques is not merely an academic exercise; it is a vital endeavor that directly expands the boundaries of what is computationally possible, leading to significant advancements in scientific understanding, engineering design, and technological innovation. In this context, the search for novel and more efficient numerical techniques, such as those leveraging machine learning, represents a key frontier in computational science.

A significant and enduring challenge in computational science is the development of robust, efficient, and unified numerical methods for approximating solutions to complex differential equations. In recent years, Artificial Neural Networks (ANNs) have emerged as a powerful tool for this purpose, demonstrating a remarkable capability to serve as universal function approximators. This inherent strength is grounded in the Universal Approximation Theorem, first established in the seminal works of Cybenko and Hornik \cite{Hornik1989,cybenko1989approximation}, which proves that a feedforward network with a single hidden layer can approximate any continuous function on a compact subset of $\mathbb{R}^n$ to any desired accuracy. The idea of leveraging ANNs as global approximators for differential equations was pioneered in the works of Lee and Kang \cite{lee1990neural}, Meade and Fernandez \cite{meade1994numerical}, Logovski \cite{logovski1992methods}, and Lagaris \cite{lagaris1998artificial}. . These foundational studies established a framework where the output of a neural network is used to construct a discrete approximation of a differential equation's solution.

The use of ANNs as numerical solvers offers several key advantages over traditional methods. The resulting network approximation is inherently smooth and continuously differentiable, a property that is highly desirable in many physical applications. Furthermore, the ANN framework operates on a set of discrete, unstructured points, which completely bypasses the complex and often time-consuming process of meshing required by conventional techniques such as the finite element method. The network output is also less susceptible to common numerical errors, including rounding, domain discretization, and approximation errors. The inherent versatility of ANNs allows the framework to be easily tuned and adapted to approximate solutions for ordinary, partial, and integral equations \cite{ngom2021fourier,lau2020oden,rao2021physics,shi2021comparative,dockhorn2019discussion}, and even for approximating Green's functions for nonlinear elliptic operators \cite{gin2021deepgreen}. These advantages are often accompanied by a reduced number of hyperparameters and the method's inherent suitability for parallel architectures.

Driven by the development of highly efficient deep learning libraries like TensorFlow and Keras \cite{terra2021keras}, a considerable body of recent literature has been published on the development of new deep learning frameworks for solving differential equations. In the present work, we contribute to this field by examining the use of a deep learning feedforward neural network for the solution of nonlinear ordinary differential equations. The framework is implemented in Python, leveraging the computational power of TensorFlow and Keras. For the purpose of illustration, two distinct boundary value problems are considered: a sixth-order linear and a nonlinear ODE. In both problems, the core methodology involves the minimization of a total loss function, which yields a highly accurate, collocated approximation of the exact solution. The network solutions are then meticulously compared with the exact solutions. 

This work builds upon the foundational research in our previous studies \cite{SMM2023,venkatachalapathy2022,venkatachalapathy2023,martinez2024approximation}, with the specific aim of addressing theoretical gaps and enhancing the computational efficiency of neural network solvers. In this investigation, we introduce and apply a highly effective mean-square metric to evaluate the performance of the proposed deep learning feedforward neural network approximation method. A key feature of our approach is the formulation of a "hard method," where the total loss function is composed of two distinct and critical components: a differential cost term that enforces the governing differential equation, and a boundary loss term that explicitly satisfies the initial and boundary conditions. This dual-term loss function ensures that the network's output is not only a valid solution but also adheres precisely to the problem's constraints.  A significant advantage of this computational implementation is that it does not require a closed-form analytical solution for the training process. Instead, the network is trained using a linear interpolation of data values to construct the tangent line and curvature, allowing the network's predictions to self-adjust to the differential equation. The overall training is further refined through backpropagation, where the network's hyperparameters are precisely fine-tuned to optimize the solution. This process yields a robust and highly accurate approximation, demonstrating the remarkable capacity of neural network frameworks to solve complex problems without the need for traditional, often limiting, numerical techniques.

The remainder of this paper is structured to systematically present our methodology and findings. 
Section~\ref{math_back} provides a comprehensive mathematical foundation, detailing the core 
method, the role of various activation functions, and a step-by-step implementation algorithm. 
In Section~\ref{ne}, we demonstrate the efficacy and robustness of the proposed 
method by applying it to two distinct numerical examples: a linear sixth-order ordinary 
differential equation (ODE) and a challenging nonlinear counterpart. Finally, 
Section~\ref{conclusions} summarizes our conclusions and discusses promising avenues for future research.

\section{Neural network architecture} \label{math_back}
A \textbf{Feedforward Neural Network (FFNN)} is employed for this study, taking an input vector $\mathbf{x} \in \mathbb{R}^n$ and producing $m$ outputs. The network architecture consists of an input layer ($l_1$), two hidden layers ($l_2$ and $l_3$), and an output layer ($l_4$), as depicted in Figure \ref{fig:NN_dia}. The input layer, $l_1$, takes the input vector $\mathbf{x} = \{x_1, x_2, \dots, x_n\}$ and connects it to the $H$ nodes of the first hidden layer, $l_2$, via a weight matrix with individual weights denoted as $w_{ij}$. The connections between the hidden layers, $l_2$ and $l_3$, are defined by weights $\tilde{w}_{ij}$, while the final hidden layer ($l_3$) is connected to the output layer ($l_4$) by weights $v_{ij}$. Each node within the network, including the hidden and output layers, receives an input that is the weighted sum of outputs from the previous layer, along with an additional term known as the  \textbf{bias}. The output of the hidden layers is then processed by an \textbf{activation function}, and the final output of the network, denoted as $N(\mathbf{x}, \mathbf{v})$, provides the approximate solution. This multi-layered structure allows the network to learn complex, non-linear relationships between the input and output data.

\begin{figure}[H]
    \centering
    \includegraphics[width=10cm]{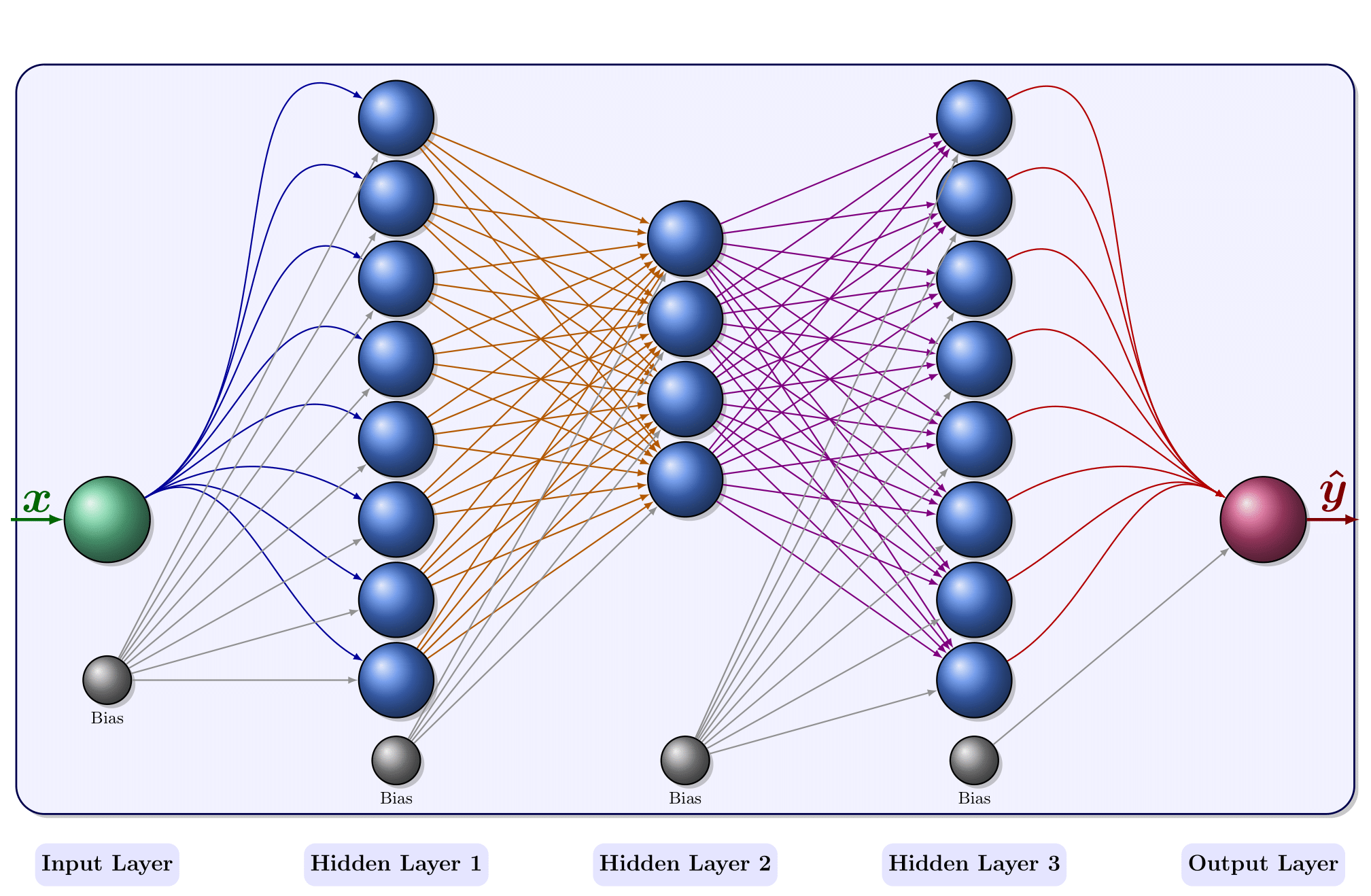}
    \caption{Feedforward neural network architecture with a single hidden layer.}
    \label{fig:NN_dia}
\end{figure}

The fundamental components that define the strength of connections within a neural network are its \textbf{weights} and \textbf{biases}, collectively referred to as the network parameters $\mathbf{P}$. The weights, in particular, determine the influence of one neuron's output on the next. A critical step in the training process is the initialization of these parameters. Initializing all weights to the same value can lead to a failure in learning, as the network's neurons would produce identical outputs and gradients, preventing the model from distinguishing different patterns. To circumvent this, the weights are initialized with random values drawn from specific distributions \cite{goodfellow2016deep}.

These weights are represented in matrix form as:
\begin{equation}
\mathbf{W} =
\begin{bmatrix}
w_{11} & w_{12} & \cdots & w_{1j}\\
w_{21} & w_{22} & \cdots & w_{2j}\\
\vdots & \vdots & \ddots & \vdots \\
w_{i1} & w_{i2} & \cdots & w_{ij}\\
\end{bmatrix}.
\end{equation}
The choice of an appropriate \textbf{initializer algorithm} is crucial, as the network's performance is highly sensitive to the initial values of its weights. Two widely used initializers are \textit{Glorot uniform} and \textit{Glorot normal} \cite{glorot}. The former draws values from a uniform distribution, while the latter uses a normal distribution. For \textit{Glorot normal}, the weights are initialized using a mean $\tilde{m}=0$ and a standard deviation $\tilde{s}=\sqrt{\frac{2}{I_{in}+ O_{out}}}$, where $I_{in}$ and $O_{out}$ are the number of input and output nodes, respectively. For \textit{Glorot uniform}, the weights are initialized from a uniform distribution with a range of $[-\sqrt{\frac{6}{I_{in}+ O_{out}}}, \sqrt{\frac{6}{I_{in}+ O_{out}}}]$.

In addition to weights, the network also incorporates a special type of constant input known as the \textbf{bias} $\mathbf{b}$. Each neuron has its own bias, which is not influenced by the outputs of other neurons. These weights and biases constitute the complete set of network parameters $\mathbf{P}$ that are iteratively adjusted during training to minimize the loss function and achieve the desired output.

\subsection{The role of activation functions}

The \textbf{activation function}, also known as a \textbf{squashing function}, is a fundamental component applied to each neuron in the hidden and output layers of the neural network. Its primary purpose is to introduce non-linearity into the network, enabling it to learn complex mappings and solve problems that cannot be addressed by a simple linear model. Mathematically, an activation function $\sigma$ is defined as a non-decreasing, real-valued function that maps from $\mathbb{R}$ to a bounded range, typically $[0,1]$ or $[-1,1]$ \cite{mcculloch1943logical}. The selection of an activation function is a critical design choice, as it significantly impacts the network's learning capacity and performance.

Among the various activation functions, this work employs several common types: \textit{Tanh}, \textit{Sigmoid}, and \textit{Softmax}. The \textbf{Sigmoid activation function}, $\sigma_s(\mathbf{x}): \mathbb{R} \to [0,1]$, is a bounded and differentiable function characterized by its distinctive "S"-shaped curve. It is defined as:
\begin{equation}
\sigma_s(\mathbf{x}) = \frac{1}{(1+ e^{-\mathbf{x}})}
\end{equation}
Its derivative, which is essential for backpropagation, is given by:
\begin{equation}
\sigma_s'(\mathbf{x}) = \sigma_s(\mathbf{x})(1- \sigma_s(\mathbf{x}))
\end{equation}

The \textbf{Softmax activation function} is particularly useful for multi-class classification problems, as it takes an $n$-dimensional input vector and produces an $n$-dimensional output vector where each value lies between 0 and 1 and the sum of all values is equal to 1. This output can be interpreted as a probability distribution over the possible classes. The function and its derivative are defined as:
\begin{align}
\sigma_t(x_i) &= \frac{e^{x_{i}}}{\sum_{j} e^{x_j}} \quad \forall i,j \in \{1, \ldots, N\} \\
\sigma_t^{\prime}(x_i) &=
\begin{cases}
\sigma_t(x_i)(1-\sigma_t(x_i)) & i = j, \\
- \sigma_t(x_i) \sigma_t(x_j) & i \neq j.
\end{cases}
\end{align}
The mathematical operation within a hidden layer involves a linear combination of the inputs and their corresponding weights, which is then passed through an activation function. The output of the $i$-th hidden node is given by:
\begin{equation}\label{eq:act_weighted}
y_{i} = \sigma(g_i) = \sigma\left(\sum_{j=1}^{n} x_{j}w_{ij}\right).
\end{equation}
This process is repeated for each hidden neuron, and the final output of the FFNN, denoted by $N_j$, is a linear combination of the activated outputs from the last hidden layer:
\begin{equation}
N_{j} = \sum_{i=1}^{H}v_{ij}y_{i}.
\end{equation}
The choice of activation function is contingent on the specific problem. The final network output is a function of the input vector and the network parameters, $N(\mathbf{x}, \mathbf{P})$, and is subsequently used to define the loss function, which serves as the metric for evaluating the solution's accuracy.

\subsection{The optimization problem and training process}

The core of training a neural network lies in solving a \textbf{minimization problem}, where the objective is to find a set of network parameters $\mathbf{P}$ that minimize the loss function $L(\mathbf{x}, \mathbf{P})$. This process is performed by \textbf{optimizers}, which are algorithms designed to iteratively update the parameters in the direction of steepest descent.

\begin{definition}[Minimization Problem]
Let $\tilde{L}: \mathbb{R}^n \to \mathbb{R}$ be a continuously differentiable function and $\mathbf{x} \in \mathbb{R}^n$ be an $n$-dimensional vector subjected to certain constraints $c$. The minimization problem for $\tilde{L}$ is then defined as follows:
\begin{equation}
\underset{\mathbf{x}\in \mathbb{R}^n} {\text{min}} \ \ \tilde{L}(\mathbf{x}) \ \ \text{subject to} \
\begin{cases}
c_i(x_1,\cdots,x_n) = r_i & i = 1,\cdots,m \\
c_j(x_1,\cdots,x_n) \geq s_j & j = 1,\cdots,z.
\end{cases}
\end{equation}
\end{definition}

\begin{definition}[Minimum]
A vector of points $\mathbf{x}^*$ is called a minimum of an objective function $\tilde{L}$ if,
\begin{equation}
\tilde{L}(\mathbf{x}^*) \leq \tilde{L}(\mathbf{x}) \ \ \forall \mathbf{x} \in \mathbb{R}^n.
\end{equation}
\end{definition}

The optimization process begins with an initial guess for the parameters and then generates a sequence of corrections until a minimum is reached. The primary role of an optimizer in a neural network is to find the optimal parameters that significantly reduce the loss function during training. A common approach is a \textbf{gradient descent-based optimization method}, which iteratively updates the parameters in the negative gradient direction of the loss function, ensuring that the parameters converge toward a minimum. The gradients are efficiently computed using the \textbf{backpropagation algorithm} \cite{backprop}, which propagates the error from the output layer back through the network, allowing for the precise adjustment of each weight and bias.

The overall workflow of the neural network model involves two distinct phases: \textbf{training} and \textbf{testing}. Initially, a training dataset is fed through the network's input layer to the hidden layers. The activation function is applied to each hidden node to produce the activated weighted sum, which is then used to define the loss function. This function serves as the primary metric for evaluating the neural network's solution quality. The training phase concludes once the optimizer has converged on a set of parameters that minimize the loss function. In the subsequent testing phase, the trained network model is evaluated on the same training data to generate the final numerical solutions, which are then used for post-processing and comparison with other results.

Here, we detail the mechanism of our neural network approach to approximate the solution to a sixth-order ODE. We consider the general form of a differential operator $\mathcal{L}$ defined on the domain $\Omega = (a, b)$:
\begin{equation}
\mathcal{L}[x,y(x), y^{\prime}(x), y^{\prime\prime}(x), \, \ldots , y^{vi}(x)] = f(x) \ \ \text{for} \ \ x \in \Omega:(a,b), \label{ode}
\end{equation}
subject to given boundary conditions at $x=a$ and $x=b$. Our approach approximates the exact solution $y(x)$ with a neural network output $\widehat{y}(\mathbf{x}, \mathbf{P})$, where $\mathbf{x}$ is the input vector and $\mathbf{P}$ represents the adjustable network parameters (weights and biases).

The primary objective of the training phase is to find the optimal parameters $\mathbf{P}$ that minimize a carefully constructed \textbf{loss function} $L(\mathbf{x},\mathbf{P})$. This loss function is a composite metric that evaluates how well the network's approximation satisfies both the differential equation and the boundary conditions. It consists of two main terms: a differential equation loss ($L_{d}$) and a boundary condition loss ($L_{bc}$).

The differential equation loss is formulated as a mean squared error, measuring the discrepancy between the network's output and the differential equation's right-hand side. This is calculated over a set of $N_{int}$ interior data points:
\begin{equation}
L_{d} = \frac{1}{N_{int}}\sum_{i=1}^{N_{int}} \left( \mathcal{L}[x_i,\widehat{y}(x_i), \widehat{y}^{\prime}(x_i), \widehat{y}^{\prime\prime}(x_i), \, \ldots, \widehat{y}^{vi}(x_i)] - f(x_i) \right)^2.
\end{equation}
Simultaneously, the boundary condition loss ensures that the network's solution adheres to the prescribed values at the boundaries of the domain. Our approach implements this as a \textbf{"hard method,"} where the boundary constraints are directly embedded into the loss function. The boundary loss is given by:
\begin{align}
L_{bc} = & (\widehat{y}(a) - y(a))^2 + (\widehat{y}^{\prime\prime}(a) - y^{\prime\prime}(a))^2 + (\widehat{y}^{iv}(a) - y^{iv}(a))^2 \notag \\
&+ (\widehat{y}(b) - y(b))^2 + (\widehat{y}^{\prime\prime}(b) - y^{\prime\prime}(b))^2 + (\widehat{y}^{iv}(b) - y^{iv}(b))^2.
\end{align}
The total loss function, which the network seeks to minimize, combines these two components into a single mean squared metric:
\begin{equation}
L = L_{d}^2+ L_{bc}^2. \label{eq:loss}
\end{equation}
This unified loss function allows the network to find an approximate solution that is not only valid across the domain but also strictly satisfies the boundary conditions, ensuring the physical consistency of the solution. 

The overall process of the neural approximation consists of a training and testing phase that can be viewed as an algorithm :
\begin{algorithm}[H]
\caption{Neural Network Solver for Boundary Value Problems}\label{alg:solver}
\begin{algorithmic}[1]
\State \textbf{Problem Formulation:} Define the sixth-order BVP in the form of a differential operator $\mathcal{L}$ as shown in Equation \eqref{ode}, along with its corresponding boundary conditions.
\State \textbf{Network Architecture Selection:} Intuitively choose a suitable feedforward neural network architecture, including the number of hidden layers, nodes, and activation functions, based on the complexity of the problem.
\State \textbf{Parameter Initialization:}
    \Statex \quad a. Initialize all network parameters $\mathbf{P} = \{\mathbf{W}, \mathbf{b}\}$ (weights and biases) using a suitable random initializer (e.g., Glorot uniform or Glorot normal).
    \Statex \quad b. Initialize a loss variable $L$ to a large value.
\State \textbf{Training Loop:} Iterate for a predefined number of epochs or until $L \le \epsilon$.
    \Statex \quad a. \textbf{Forward Pass:} Propagate the input data points through the network to compute the network's output $\widehat{y}(\mathbf{x}, \mathbf{P})$.
    \Statex \quad b. \textbf{Loss Computation:} Calculate the total mean squared loss $L = L_d^2 + L_{bc}^2$, as given by Equation \eqref{eq:loss}.
    \Statex \quad c. \textbf{Backpropagation:} Compute the gradients of the loss function with respect to each network parameter ($\frac{\partial L}{\partial \mathbf{P}}$).
    \Statex \quad d. \textbf{Parameter Update:} Update the network parameters using an optimizer: $\mathbf{P}_{n+1} = \mathbf{P}_n - \eta \frac{\partial L}{\partial \mathbf{P}_n}$, where $\eta$ is the learning rate.
\State \textbf{Convergence Check:} If $L \le \epsilon$, or the maximum number of epochs has been reached, terminate the training process.
\State \textbf{Testing and Evaluation:}
    \Statex \quad a. Feed the testing dataset into the trained network with the saved, optimized parameters.
    \Statex \quad b. The network's output on this dataset constitutes the numerical solution for the BVP.
    \Statex \quad c. Perform post-processing and analysis by comparing the obtained solutions with reference data from other numerical or analytical methods.
\end{algorithmic}
\end{algorithm}
\section{Numerical Experiments}

This section presents the application of our proposed feedforward neural network method to two distinct boundary value problems: sixth-order linear and nonlinear ODEs. A key objective of these experiments is to demonstrate the point-wise convergence and high accuracy of the neural network's approximate solution. Since both problems lack a known analytical solution, the performance of our neural network solver is rigorously validated by comparing its output against exact solutions. The entire training, testing, and optimization procedure for the network was conducted on personal laptops.  This computational power allowed for the efficient execution of complex network architectures and extensive training over more than 30,000 epochs, with each run completing in a matter of minutes. For all numerical experiments, a consistent learning rate of $0.001$ was utilized. Extensive testing was performed to identify the optimal combination of activation functions and optimizers for each problem. The subsequent subsections detail the specific results obtained for the two test cases, providing a comprehensive evaluation of the framework's performance.\\

The core of our solver is a fully connected, feedforward neural network, commonly known as a \textit{multi-layer perceptron}. The architecture and training hyperparameters are specifically configured as follows:

\begin{itemize}
    \item \textbf{Network Architecture:} The network consists of an input layer with a single neuron for the spatial coordinate $x$, one hidden layer containing 16 neurons to approximate the complexity of the solution, and an output layer with a single neuron that provides the final solution, $\widehat{y}(x)$.

    \item \textbf{Activation Functions:} The hyperbolic tangent ($\tanh$) function is utilized as the activation for the hidden layer, as its smoothness is beneficial for approximating differential equations. The output neuron employs a linear activation function to allow for an unbounded solution space.

    \item \textbf{Initialization Strategy:} Network weights are initialized using the Glorot Normal distribution (also known as Xavier initialization), a standard practice for deep networks that helps prevent vanishing or exploding gradients. All bias terms are initialized to zero.

    \item \textbf{Optimization Algorithm:} The Adamax optimization algorithm is used to update the network parameters during training, with a constant learning rate set to $1 \times 10^{-3}$.

    \item \textbf{Physics-Informed Loss Function:} A custom loss function guides the training process by encoding the physics of the problem. This function is the sum of two components: the mean squared error of the ODE residual, and the mean squared error from the specified boundary conditions at $x=0$ and $x=1$.

    \item \textbf{Training Procedure:} The model is trained for a total of 13,000 epochs to iteratively minimize the composite loss function until a satisfactory level of convergence is achieved.
\end{itemize}

\noindent \textbf{Example-1: }
In this example, we apply the proposed neural solver to approximate the solution for the following linear, sixth-order boundary value problem.

\begin{subequations}
\begin{align}
    & y^{(6)}(x) - y(x) = -6e^x, \quad \text{for} \quad 0 < x < 1, \\
    & y(0) = 0, \quad y^{(2)}(0) = -1, \quad y^{(4)}(0) = -3, \\
    & y(1) = 0, \quad y^{(2)}(1) = -2e, \quad y^{(4)}(1) = -4e.
\end{align}
\end{subequations}
The exact solution to the BVP is given by $y(x) = (1 - x)e^x$. The composite loss function $\mathcal{L}$ used to train the neural network approximation, $\widehat{y}(x)$, is defined as the sum of the differential equation loss ($\mathcal{L}_{d}$) and the boundary condition loss ($\mathcal{L}_{b}$).

\begin{align*}
    \text{Differential Equation Loss: } \quad \mathcal{L}_{d} &= \frac{1}{N_{\text{int}}}\sum_{i=1}^{N_{\text{int}}} \left( \widehat{y}^{(6)}(x_i) - \widehat{y}(x_i) + 6e^{x_i} \right)^2 \\
    \\
    \text{Boundary Loss: } \quad \mathcal{L}_{b} &= \left( \widehat{y}(0) \right)^2 + \left( \widehat{y}^{(2)}(0) + 1 \right)^2 + \left( \widehat{y}^{(4)}(0) + 3 \right)^2 \\
    &\quad + \left( \widehat{y}(1) \right)^2 + \left( \widehat{y}^{(2)}(1) + 2e \right)^2 + \left( \widehat{y}^{(4)}(1) + 4e \right)^2 \\
    \\
    \text{Total Loss: } \quad \mathcal{L} &= \mathcal{L}_{d} + \mathcal{L}_{b}.
\end{align*}

The subsequent table presents a detailed comparison between the solution approximated by our neural network and the exact analytical solution. An analysis of the results reveals that our neural network's predictions demonstrate a higher degree of accuracy than the numerical solutions reported in \cite{wazwaz2001numerical}. This underscores the enhanced performance and precision of our proposed method for this particular problem.

\begin{table}[H]
\centering
\caption{Comparison of analytical and numerical solutions with absolute error.}
\label{tab:solution_comparison}
\begin{tabular}{cccc}
\toprule
\textbf{x} & \textbf{Analytical} & \textbf{Numerical} & \textbf{Error} \\
\midrule
0.0 & 1.00000000 & 1.00000000 & 0.00000000 \\
0.1 & 0.99465383 & 0.99458730 & 0.00006652 \\
0.2 & 0.97712221 & 0.97701275 & 0.00010945 \\
0.3 & 0.94490117 & 0.94476765 & 0.00013351 \\
0.4 & 0.89509482 & 0.89495212 & 0.00014270 \\
0.5 & 0.82436064 & 0.82422036 & 0.00014028 \\
0.6 & 0.72884752 & 0.72871995 & 0.00012757 \\
0.7 & 0.60412581 & 0.60401964 & 0.00010617 \\
0.8 & 0.44510819 & 0.44503129 & 0.00007690 \\
0.9 & 0.24596031 & 0.24591930 & 0.00004101 \\
1.0 & 0.00000000 & 0.00000000 & 0.00000000 \\
\bottomrule
\end{tabular}
\end{table}

\begin{figure}[H]
    \centering
    \includegraphics[width=15cm]{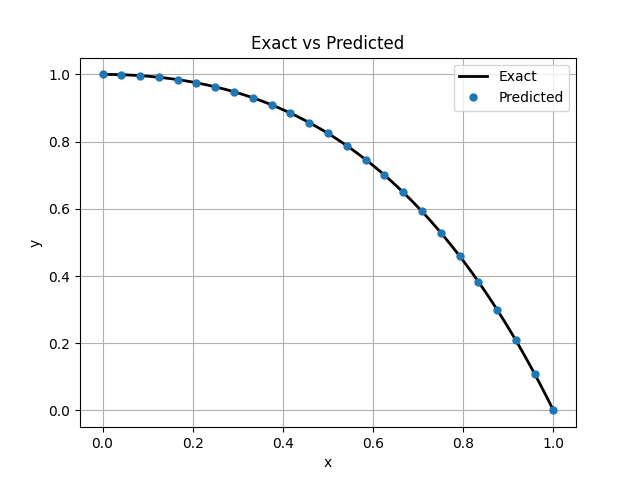}
\caption{A graphical comparison between the exact solution (solid line) and the solution approximated by the neural network (markers). The close agreement between the two illustrates the high accuracy of our proposed model.}
    \label{fig:NN_dia}
\end{figure}

\noindent \textbf{Example-2: }
In this example, we apply the proposed neural solver to approximate the solution for the following linear, sixth-order boundary value problem.

\begin{subequations}
\begin{align}
    & y^{(6)}(x) - e^{-x} y^2(x) =0, \quad \text{for} \quad 0 < x < 1, \\
    & y(0) = 1, \quad y^{(2)}(0) = 1, \quad y^{(4)}(0) = 1, \\
    & y(1) = e, \quad y^{(2)}(1) = e, \quad y^{(4)}(1) = e.
\end{align}
\end{subequations}
The exact solution to the BVP is given by $y(x) = e^x$. The composite loss function $\mathcal{L}$ used to train the neural network approximation, $\widehat{y}(x)$, is defined as the sum of the differential equation loss ($\mathcal{L}_{d}$) and the boundary condition loss ($\mathcal{L}_{b}$).

\begin{align*}
    \text{Differential Equation Loss: } \quad \mathcal{L}_{d} &= \frac{1}{N_{\text{int}}}\sum_{i=1}^{N_{\text{int}}} \left( \widehat{y}^{(6)}(x_i) -  e^{-x_i }\widehat{y}^{2}(x_i)  \right)^2 \\
    \\
    \text{Boundary Loss: } \quad \mathcal{L}_{b} &= \left( \widehat{y}(0)  -1 \right)^2 + \left( \widehat{y}^{(2)}(0) - 1 \right)^2 + \left( \widehat{y}^{(4)}(0) -  \right)^2 \\
    &\quad + \left( \widehat{y}(1) - e\right)^2 + \left( \widehat{y}^{(2)}(1) - e \right)^2 + \left( \widehat{y}^{(4)}(1) - e \right)^2 \\
    \\
    \text{Total Loss: } \quad \mathcal{L} &= \mathcal{L}_{d} + \mathcal{L}_{b}.
\end{align*}

The subsequent table presents a detailed comparison between the solution approximated by our neural network and the exact analytical solution. An analysis of the results reveals that our neural network's predictions demonstrate a higher degree of accuracy than the numerical solutions reported in \cite{wazwaz2001numerical}. This underscores the enhanced performance and precision of our proposed method for this particular problem.

\begin{table}[H]
\centering
\caption{Comparison of the analytical solution, the Artificial Neural Network (ANN) solution, and the absolute error.}
\label{tab:ann_solution_comparison}
\begin{tabular}{cccc}
\toprule
\textbf{x} & \textbf{Analytical Solution} & \textbf{ANN Solution} & \textbf{Error} \\
\midrule
0.0 & 1.00000000 & 1.00821758 & 0.00821758 \\
0.1 & 1.10517092 & 1.11291075 & 0.00773983 \\
0.2 & 1.22140276 & 1.22851976 & 0.00711700 \\
0.3 & 1.34985881 & 1.35623278 & 0.00637398 \\
0.4 & 1.49182470 & 1.49736021 & 0.00553551 \\
0.5 & 1.64872127 & 1.65334737 & 0.00462610 \\
0.6 & 1.82211880 & 1.82578868 & 0.00366988 \\
0.7 & 2.01375271 & 2.01644319 & 0.00269048 \\
0.8 & 2.22554093 & 2.22725196 & 0.00171103 \\
0.9 & 2.45960311 & 2.46035723 & 0.00075411 \\
1.0 & 2.71828183 & 2.71812365 & 0.00015818 \\
\bottomrule
\end{tabular}
\end{table}

\begin{figure}[H]
    \centering
    \includegraphics[width=15cm]{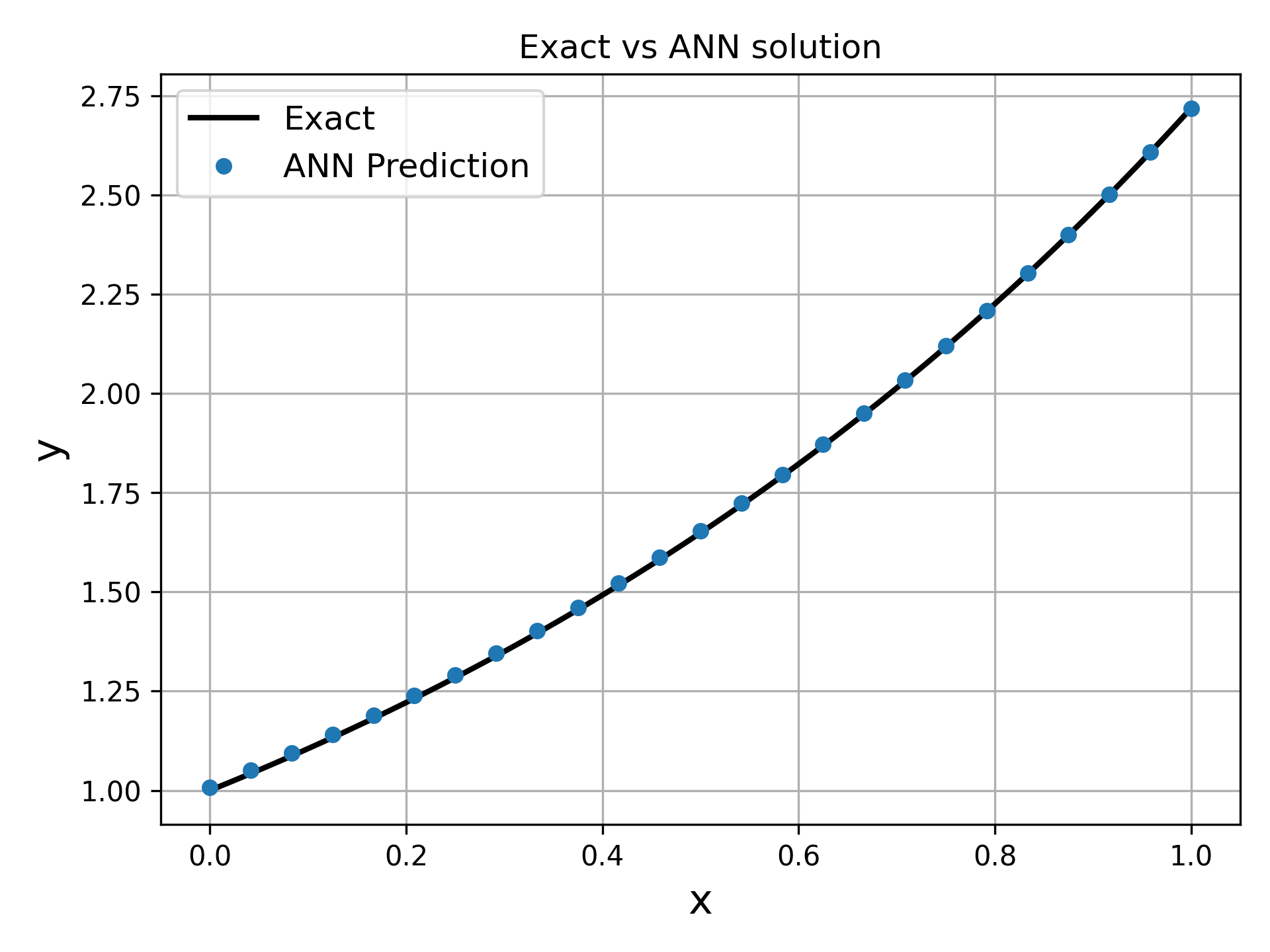}
\caption{A graphical comparison between the exact solution (solid line) and the solution approximated by the neural network (markers). The close agreement between the two illustrates the high accuracy of our proposed model.}
    \label{fig:NN_dia}
\end{figure}


\section{Conclusions} \label{conclusions}
In this paper, we introduced a deep neural network framework as a powerful tool for approximating the solutions of high-order ordinary differential equations. Our approach is built on a straightforward and easy-to-implement physics-informed methodology. A single mean-squared error loss function effectively unifies the differential equation's structure and its boundary conditions, which simplifies the problem-solving process. We also conducted a systematic investigation into the effects of various hyperparameters, such as activation functions, optimizers, and network architecture, to optimize the model's performance. The efficacy of our framework was demonstrated on two challenging sixth-order boundary value problems: one linear and one nonlinear. In both cases, the neural network yielded solutions with a high degree of accuracy, closely matching the exact or established benchmark solutions. These results validate our method's capability to handle the complexity and stiffness often associated with high-order differential equations without requiring problem-specific modifications.

Our findings suggest that this neural network approach is a compelling and robust alternative to classical numerical solvers like the Runge-Kutta or finite difference methods. A key advantage of our method is its simplicity; it circumvents the need for complex grid generation, iterative adjustments for boundary conditions, or linearization techniques for nonlinear problems. The solution is found directly by minimizing a well-defined objective function.

We plan to expand this framework to tackle various complex problems, such as systems involving fractional-order derivatives and partial differential equations. Additionally, establishing a more rigorous theoretical foundation for this approach is a crucial next step. Future investigations will focus on proving the convergence of the ANN method and formally defining the relationship between the loss function and the true solution error.

\section*{Acknowledgements}
The author is deeply grateful to Dr. S. M. Mallikarjunaiah (Associate Professor of Mathematics, Texas A\&M University–Corpus Christi) for proposing this research problem and for his invaluable guidance throughout the development of the methodology and neural network code. The valuable assistance of Ms. Pavithra Venkatachalapathy (PhD student, Texas Tech University) is also gratefully acknowledged.


\bibliographystyle{model3-num-names}
\bibliography{ML_references}

\end{document}